\documentclass[12pt]{article}
\usepackage{graphicx,amssymb,amsmath,amsthm}
\usepackage{comment,cite}

\newtheorem{theorem}{Theorem}

\newtheorem{definition}{Definition}

\newtheorem{property}{Property}
\newtheorem{remark}{Remark}
\newtheorem{corollary}{Corollary}

\def\N{\mathbb{N}}

\DeclareMathOperator{\Span}{span}
\DeclareMathOperator{\supp}{supp}
\def\Del{\Delta}
\def\tDel{\tilde{\Del}}
\def\D{{\cal D}}
\def\P{{\cal P}}
\newcommand{\JI}{J_{\rm I}}
\newcommand{\JL}{J_{\rm L}}
\newcommand{\JR}{J_{\rm R}}
\newcommand{\Pb}[1][b]{\P_{#1}(c-mb,d)}
\newcommand{\PLb}[1][b]{\P_{#1}(c-mb,c)}
\newcommand{\PIb}[1][b]{\P_{#1}[c,d-mb]}
\newcommand{\PRb}[1][b]{\P_{#1}(d-mb,d)}
\newcommand{\Pbp}[1][b']{\P_{#1}(c-mb',d)}

\newcommand{\JDb}{\Pb}
\newcommand{\JIDb}{\PIb}
\newcommand{\JLDb}{\PLb}
\newcommand{\JRDb}{\PRb}

\def\be{\begin{equation}}
\def\ee{\end{equation}}
\def\la{\langle}
\def\ra{\rangle}
\def\alt{\tilde\alpha}
\title{Cardinal B-spline dictionaries on a compact interval}
\author{Miroslav Andrle and Laura Rebollo-Neira\thanks{Support from EPSRC (GR$/$R86355$/$01) is acknowledged.}\\
NCRG, Aston University, Birmingham B4 7ET, UK}
\date{}
\begin{document}
\maketitle
\begin{abstract}
A prescription for constructing dictionaries for cardinal
spline spaces on a compact interval is provided.
It is proved that such
spaces can be spanned by dictionaries which are
built by translating
a prototype B-spline function of fixed support into
the knots of the required cardinal spline space.
This implies that cardinal spline spaces on a
compact interval can be spanned by dictionaries of
cardinal B-spline functions of broader support that
the corresponding basis function.  
\end{abstract}

\section{Introduction}
The problem of signal approximation outside
the orthogonal basis setting is a non-linear
problem, per se, in the following sense:
Let us consider that
a signal $f$ in a separable Hilbert space,
equipped with an inner product
$\la \cdot, \cdot \ra$
and a norm $||f||= \la f,f \ra^{1/2}$, is to
be approximated as the
$M$-term superposition
\be
f^M= \sum_{n=1}^M c_n^M \alpha_n,
\label{at}
\ee
where 
$\{\alpha_n\}_{n=1}^{M}$ are fixed
elements of a non-orthogonal basis, which  if   
well localized are often refereed to as 
{\em{atoms}}. In order to
construct the approximation $f^M$ of $f$ minimizing the
distance
$||f-f^M||$ 
the coefficients $c_n^M$ in \eqref{at} should be computed as
$c_n= \la \alt^M_n ,f\ra$, where the dual sequence
$\{\alt_n^M\}_{n=1}^{M}$ is biorthogonal to
$\{\alpha_n\}_{n=1}^{M}$, and, in addition, the superscript
is meant to indicate that 
$\Span \{\alpha_n\}_{n=1}^{M} \equiv \Span
\{\alt_n^M\}_{n=1}^{M}$. Hence,
in order to account for
the inclusion (or respective elimination) of one
term in \eqref{at}, all the elements of the
dual sequence need to be modified for the  coefficients 
 of the  new  
approximation to minimize the distance to the signal $f$.
Such modifications can be performed in an effective 
manner by means of
adaptive biorthogonalization techniques \cite{fobio,babio}. 
However, the problem of choosing the $M$ elements of a
non-orthogonal basis best representing the signal is 
a highly nonlinear one \cite{dav,dev,ti}. 
Since the general
problem has not feasible solution in polynomial time,
in some practical situations it is addressed by algorithms 
which evolve by fixing atoms at each iteration step. These 
 approaches are known as adaptive pursuit strategies
\cite{ma,dav,pati,oomp,boomp}. 
They operate by selecting atoms 
from a redundant set, called a
{\em{dictionary}}\cite{ma}.  

It should be stressed that, as far as complexity
is concerned, 
there is not much to lose by applying the above mentioned
 practical selection strategies on a redundant
 dictionary, rather than 
on a non-orthogonal basis. Nevertheless, 
as has been shown in the context 
 of several applications, there is potentially much to win
 in relation to sparseness 
 of the signal representation. Hence the
 motivation of this Communication: 
{\em {We prove here that 
 cardinal spline spaces on a compact interval are amenable to
 be spanned by dictionaries of cardinal B-spline functions
 of broader support than the corresponding basis functions.}}

It goes without saying that 
splines have been used with success 
in wavelet theory 
and applications to signal processing
\cite{chui1,chui, chui3, unser}. 
In particular, the construction of multiresolution
based spline wavelets on a bounded interval
is explained in great detail in 
\cite{chui, chui2, quak}. 
Here we focus on B-spline functions 
and discuss the way of 
 going from  B-spline 
basis on the interval to B-spline dictionaries 
which are endowed with a very
 interesting 
property. 
Cardinal spline spaces on a bounded interval are finite
dimensional linear spaces. 
The dimension is actually
given by the number $d=m+N$, where $m$ is the order 
of the splines being considered and $N$  the number of
 knots in the interval. By fixing the order 
of splines, the usual way of increasing the subspace 
dimension is to increase the 
number of knots, decreasing thereby the distance,
$b$, say, between two adjacent knots. Let us recall that
the basis functions for cardinal spline spaces
(B-spline basis) 
have compact support which, except for the boundary
functions, is of length  $mb$. Thus, 
by decreasing the distance $b$ between  knots, 
the support of the functions is reduced.   
With the exception of the boundary functions 
\cite{chui, chui2}
a B-spline basis for the subspace arises by successive 
translations of a
prototype function. 
The translation parameter is the distance between knots.
It is clear then that if 
the translation operation is carried out with a
translation 
parameter $b'$, such that $b/b'$ is an integer,
some redundancy will be introduced. 
However, the main result of this contribution is to
prove that, by such a
procedure, one can generate the span for the cardinal 
spline space  associated with the distance $b'$ 
between knots. 
We believe this to be a
remarkable feature of splines, which is potentially useful
in relation to signal representation.
Such a possibility
will be illustrated using matching pursuit strategies for 
representing a signal by selecting atoms
from our cardinal spline dictionaries. 

The paper is organized as follows: 
Section~\ref{sec:splines} provides some definitions 
and background on splines
relevant for our purpose.
Section~\ref{sec:dictionary} gives the proof of the main
theorem establishing the above mentioned
property of the proposed dictionaries.
The potential suitability of such dictionaries 
for recursive signal approximation is illustrated in 
Section~\ref{sec:application}. The conclusions are drawn in 
Section~\ref{sec:conclusion}.

\section{Splines on a compact  interval}
\label{sec:splines}
We introduce here some notation and basic definitions 
which are relevant for our purpose. 
For an in depth treatment 
on splines  we refer to \cite{boor,schum}. 
\begin{definition}
Given a compact interval $[c,d]$ we define 
a {\em partition} of $[c,d]$ as the finite set of points
\begin{equation}\label{Delta}
\Del=\{x_i\}_{i=0}^{N+1}, N\in\N,\,\,\text{such that}
\,\,
c=x_0<x_1<\cdots<x_{N}<x_{N+1}=d.
\end{equation}
We further define $N$ subintervals $I_i, i=0,\dots,N$ as: 
$I_i=[x_i,x_{i+1}), i=0,\dots,N-1$ and $I_N=[x_N,x_{N+1}]$. 
\end{definition}
\begin{definition}\label{splinespace}
Let $\Pi_{m}$ be the 
space of polynomials  of
degree smaller or equal to $m\in\N_0=\N\cup\{0\}$. 
Let $m$ be a positive integer and define
\begin{equation}
S_m(\Del)=\{f\in C^{m-2}[c,d]\ ; \ f|_{I_i}\in\Pi_{m-1}, i=0,\dots,N\},
\end{equation}
where  $f|_{I_i}$ indicates the
restriction of the function $f$ on the
interval ${I_i}$.

We call $S_m(\Del)$ the 
{\em space of polynomial splines (or splines) 
of order $m$ with simple knots
at the points $x_1,\dots,x_N$.} 
\end{definition}
Let us recall two well known properties of $S_m(\Del)$.
\begin{property}
$S_m(\Del)$ is a linear space of dimension $m+N$ 
\cite[Theorem 4.4]{schum}.
\end{property}
Moreover, it readily follows 
from Definition \ref{splinespace} that
\begin{property}\label{prop1}
If $\Del$ and $\Del'$ 
are two partitions of the interval $[c,d]$
such that $\Del\subset\Del' $, then  
 $S_m(\Del)\subset S_m(\Del'), m\in \N.$
\end{property}
In order to construct a particular basis for $S_m(\Del)$ 
it is necessary to introduce the so-called 
{\em extended partition}.
\begin{definition}
Let $\Del$ be a partition 
of $[c,d]$  and 
let us consider
\begin{equation*}
y_1\leq y_2\leq\cdots \leq y_{2m+N}
\end{equation*}
 such that
\begin{equation}
y_1\leq \cdots \leq y_{m} \leq c,\quad
d \leq y_{m+N+1}\leq \cdots \leq y_{2m+N} 
\end{equation}
and
\begin{equation*}
y_{m+1}<\cdots <y_{m+N}=x_1,\dots,x_N.
\end{equation*}
We call $\tDel=\{y_i\}_{i=1}^{2m+N}$ an
{\em extended partition with single inner knots associated with}
$S_m(\Del)$.
\end{definition}
The points $\{y_i\}_{i=m+1}^{m+N}$ in an extended
partition $\tDel$ associated with $S_m(\Del)$ are 
uniquely determined, 
however the first and last $m$ points in $\tDel$ 
can be chosen arbitrarily.

With each fixed extended partition $\tDel$ there is associated a
unique B-spline basis for $S_m(\Del)$, that we denote as
$\{B_i\}_{i=1}^{m+N}$. Full details on how to 
construct such basis are given in \cite[Theorem 4.9]{schum} and 
\cite{chui}. The basis functions $B_i$ satisfy for $x\in [c,d]$
\begin{align}
B_i(x)&=0\qquad &&\text{if } x\notin [y_i,y_{i+m}], \\
 B_i(x)&>0 \qquad &&\text{if } x\in (y_i,y_{i+m}).\label{support}
\end{align}
In the case of equally spaced knots the corresponding splines are called cardinal. Moreover all the cardinal 
B-splines of order $m$ can be obtained from 
one cardinal B-spline $B(x)$ associated with the uniform simple knot sequence
 $0,1,\dots,m$. 
Such a function is given as
\begin{equation}
B(x)=\frac{1}{m!}\sum_{i=0}^m(-1)^i\binom{m}{i}(x-i)^{m-1}_+,
\end{equation}
where $(x-i)^{m-1}_+$ is equal to $(x-i)^{m-1}$ if $x-i>0$ and 0 otherwise.
If $y_i,\dots,y_{i+m}$ are equally spaced,  
then 
\begin{equation}
B_i(x)=\frac{1}{b}B\left(\frac{x-y_i}{b}\right),
\end{equation}
where $b$ is the distance between two adjacent knots.

Let as recall that, given a partition $\Del$ and a
 real number $r$  by the operation $\Del+r$
we obtain:
\begin{equation*}
\Del+r=\{x_i+r\}_{i=0}^{N+1}=
\{x_0+r,x_1+r,\dots,x_{N+1}+r\}.
\end{equation*}
In order to retain the boundary points $x_0$ and
$x_{N+1}$ we define a new operation `$\Del \uplus r$' as
follows:
\begin{definition} Given a partition $\Del$ and
a real number $r$, $0<r<\min_i(x_{i+1}-x_i), \ i=0,\dots,N$ by the operation $\Del\uplus r$
we obtain:
\begin{equation*}
\Del\uplus r= \{x_i+r\}_{i=0}^{N}  
\cup \{x_0,x_{N+1}\}=
\{x_0,x_0+r,\dots,x_N+r,x_{N+1}\}.
\end{equation*}
\end{definition}
Since we are interested only 
in cardinal spline spaces 
hereafter we will 
consider an equidistant partition
of $[c,d]$. Such a partition is thereby uniquely 
determined by the interval $[c,d]$ and
the distance, $b$ say, between two adjacent points. 
It is assumed that the interval $[c,d]$ contains at least 
 one complete B-spline function, i.e., $d-c\geq mb$.
The definition of
equidistant partition is extended to open/semi-open intervals 
as indicated below.
\begin{definition}
We construct an 
equidistant partition of  $[c,d], (c,d), (c,d],$ and 
$[c,d)$ with distance $b$ between adjacent points 
(such that $(d-c)/b=$integer)
as follows:
\begin{align}
\P_b[c,d]&= \{c,c+b,\dots,d-b,d\}, \qquad && 
\P_b(c,d)=\{c+b,\dots,d-b\}, \nonumber\\
\P_b(c,d]&=\{c+b,\dots,d-b,d\}, \qquad &&
\P_b[c,d)=\{c,c+b,\dots,d-b\} \nonumber.
\end{align}
\end{definition} 
As already mentioned, the selection of a particular 
extended partition $\tDel$  yields a particular 
B-spline basis for  $S_m(\Del)$. Two  possible  
choices for the first and last $m$ points determining the 
extended partition are the following:
\begin{itemize}
\item[ i)]The points $\{y_i\}_{i=1}^{m} $  
and $\{y_i\}_{i=N+m+1}^{2m+N} $ are determined 
in order for the whole sequence 
$\tDel=\{y_i\}_{i=1}^{2m+N}$ to be endowed with 
the equidistant property, i.e., $\tDel=\P_b(c-mb,d+mb)$.
Such a sequence is called an  
{\em equally spaced extended partition} (ESEP). 
\item[ii)]The points lying outside $[c,d]$ are given the value
of the closest point in the interval, 
i.e., 
$y_1= \cdots = y_m = c$ and $y_{m+N+1}= \cdots = y_{2m+N}= d$.
This is called an 
{\em extended partition with $m$-tuple knots on the border} 
(EPKB). 
\end{itemize}
\section{Building B-spline dictionaries}
\label{sec:dictionary}
\begin{figure}[!ht]
\begin{center}
\includegraphics[width=6cm]{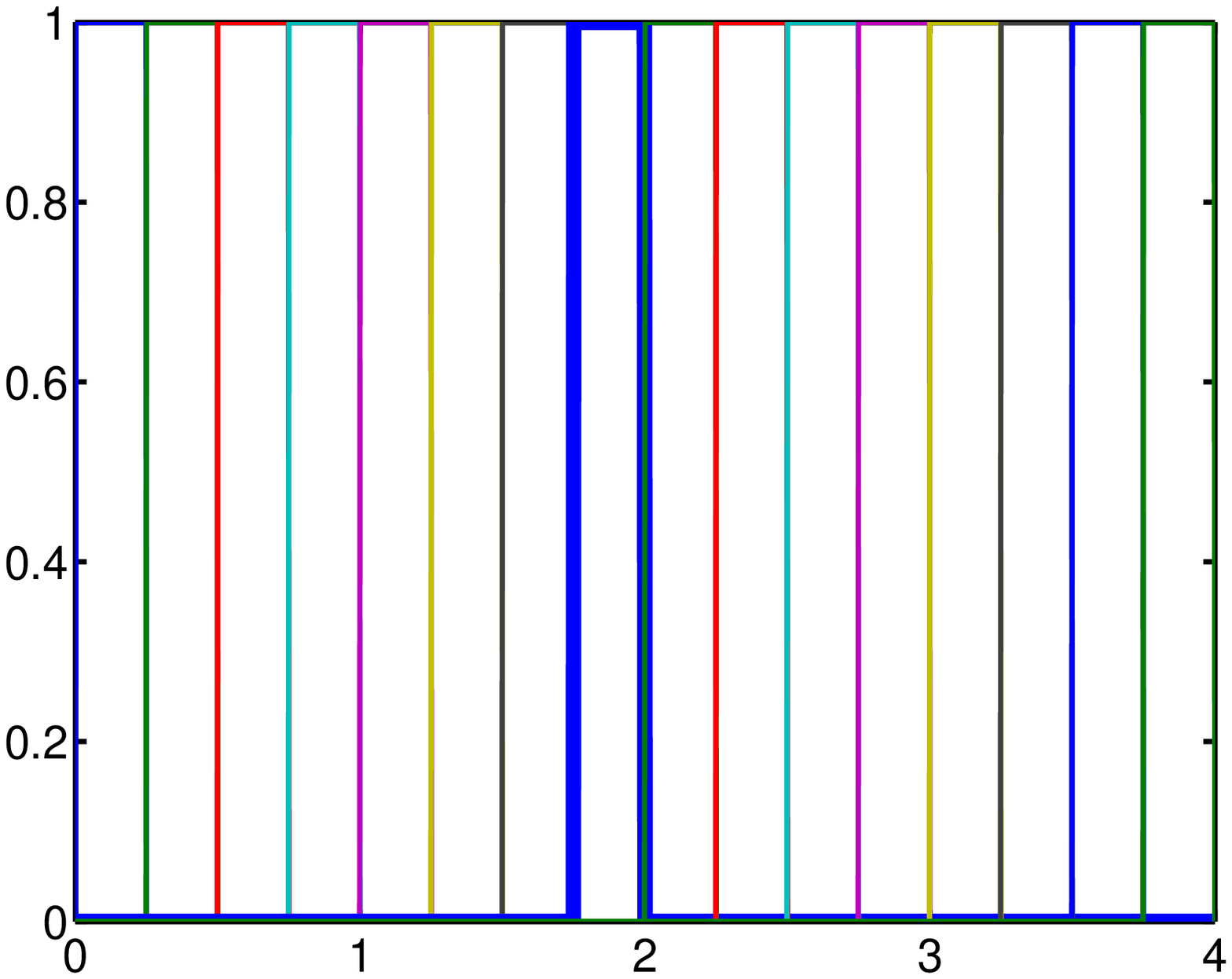}\qquad
\includegraphics[width=6cm]{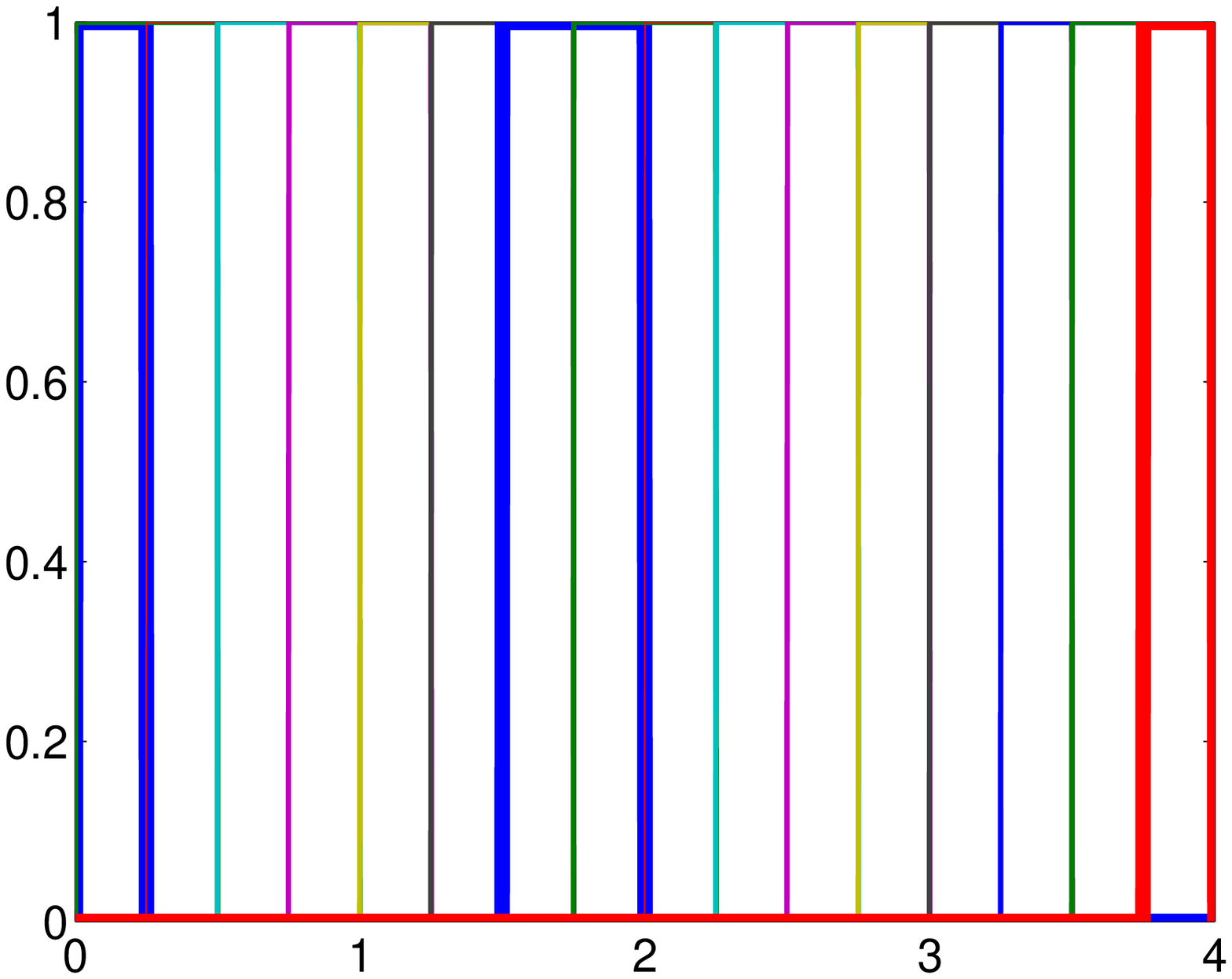}
\includegraphics[width=6cm]{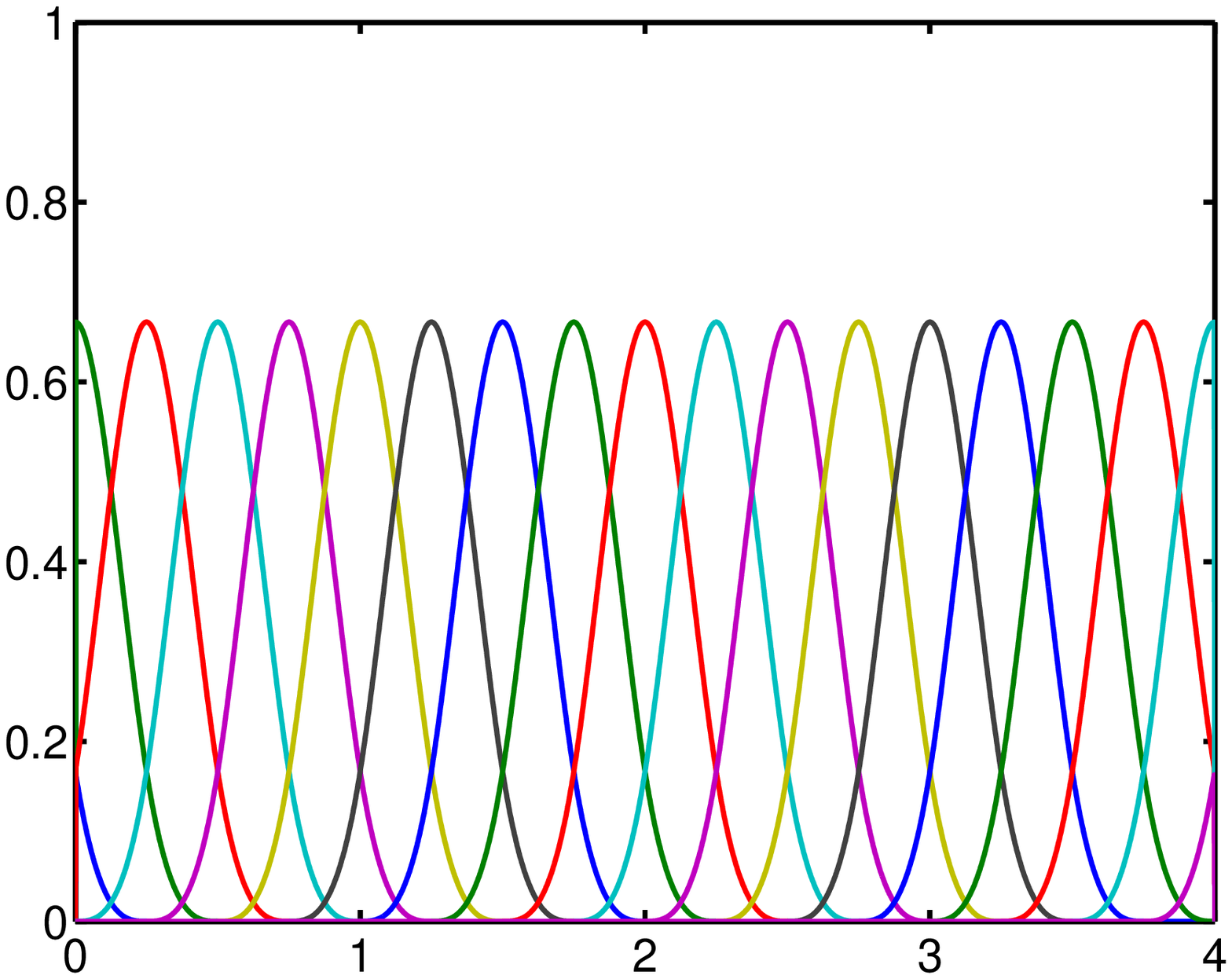}\qquad
\includegraphics[width=6cm]{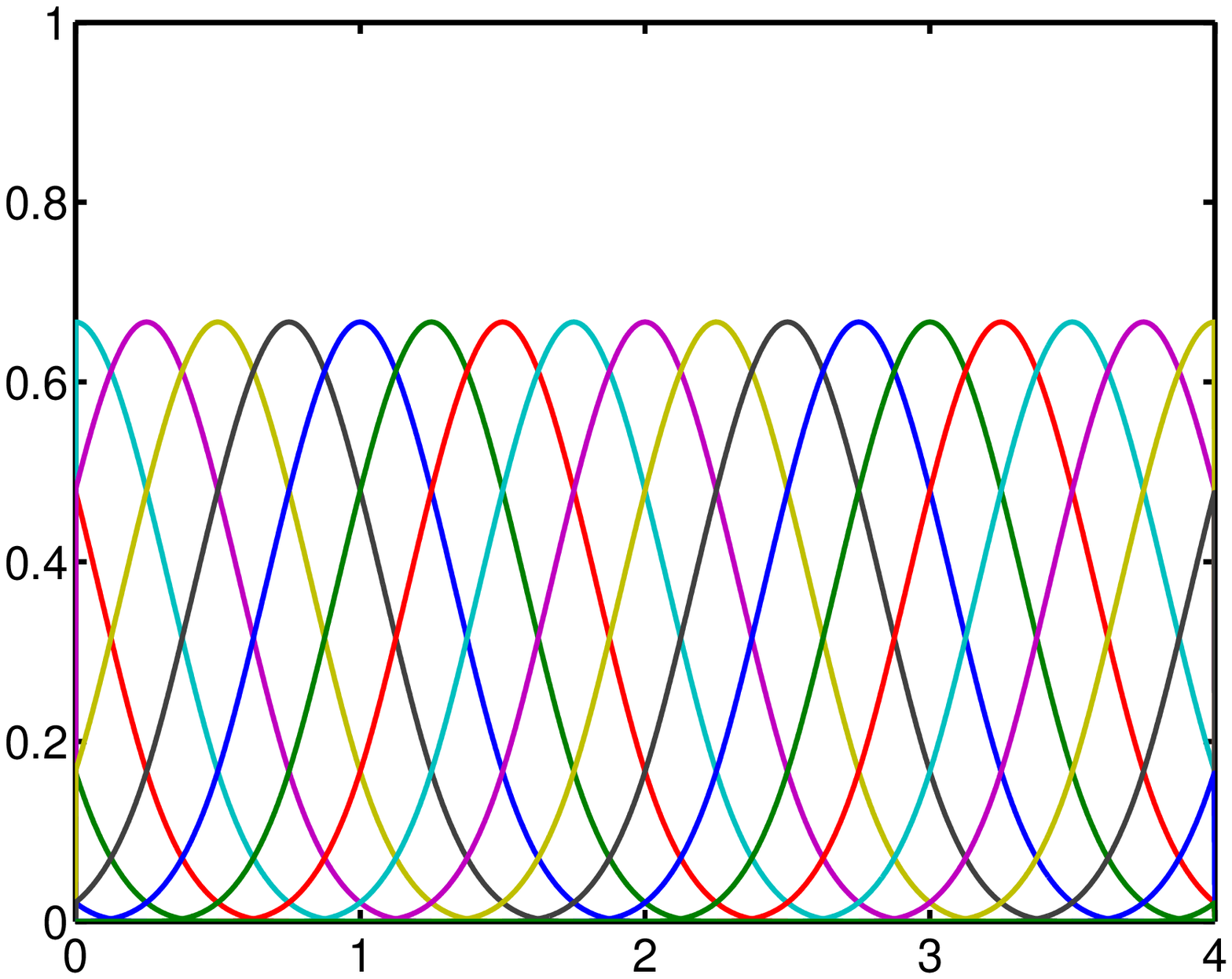}
\includegraphics[width=6cm]{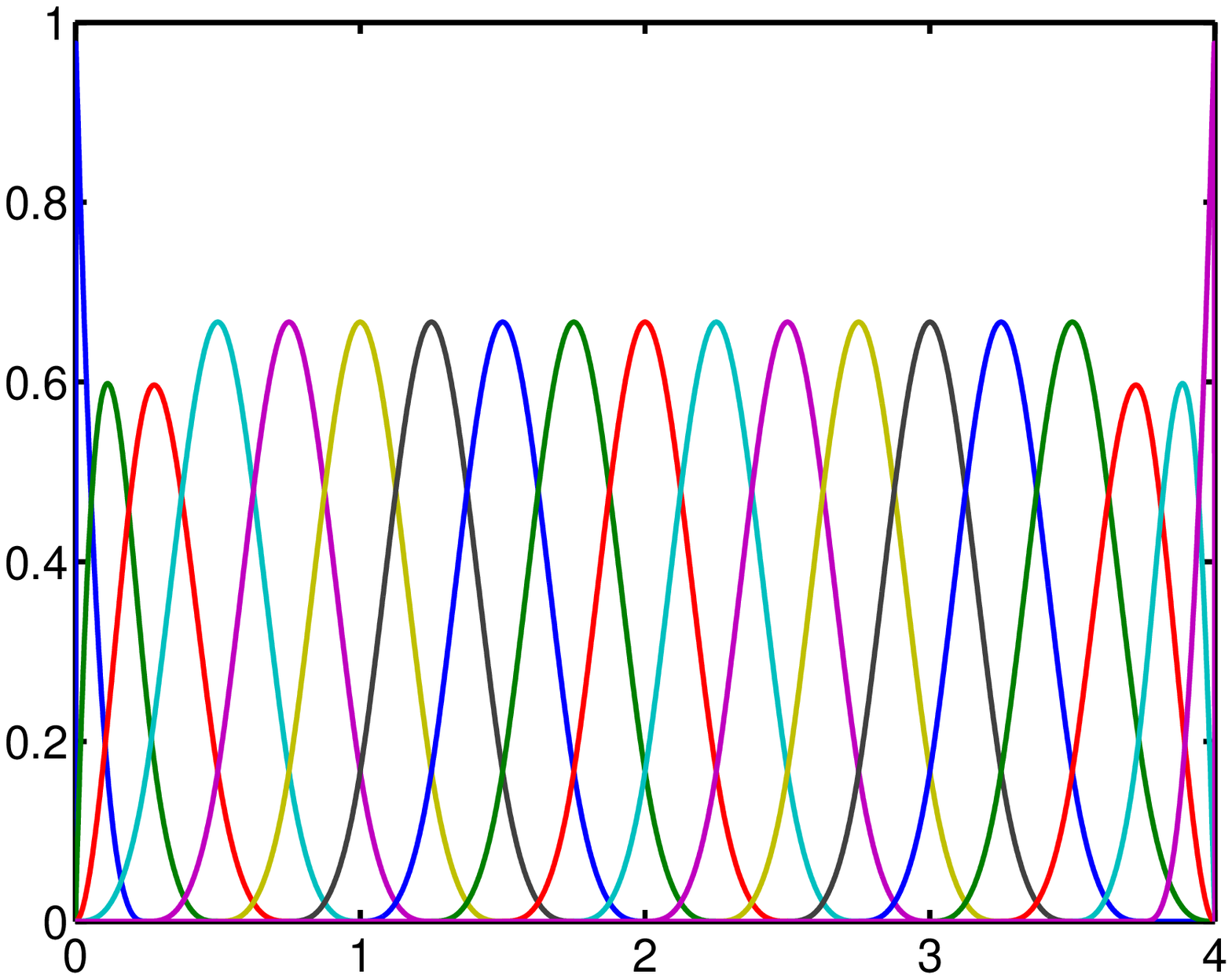}\qquad
\includegraphics[width=6cm]{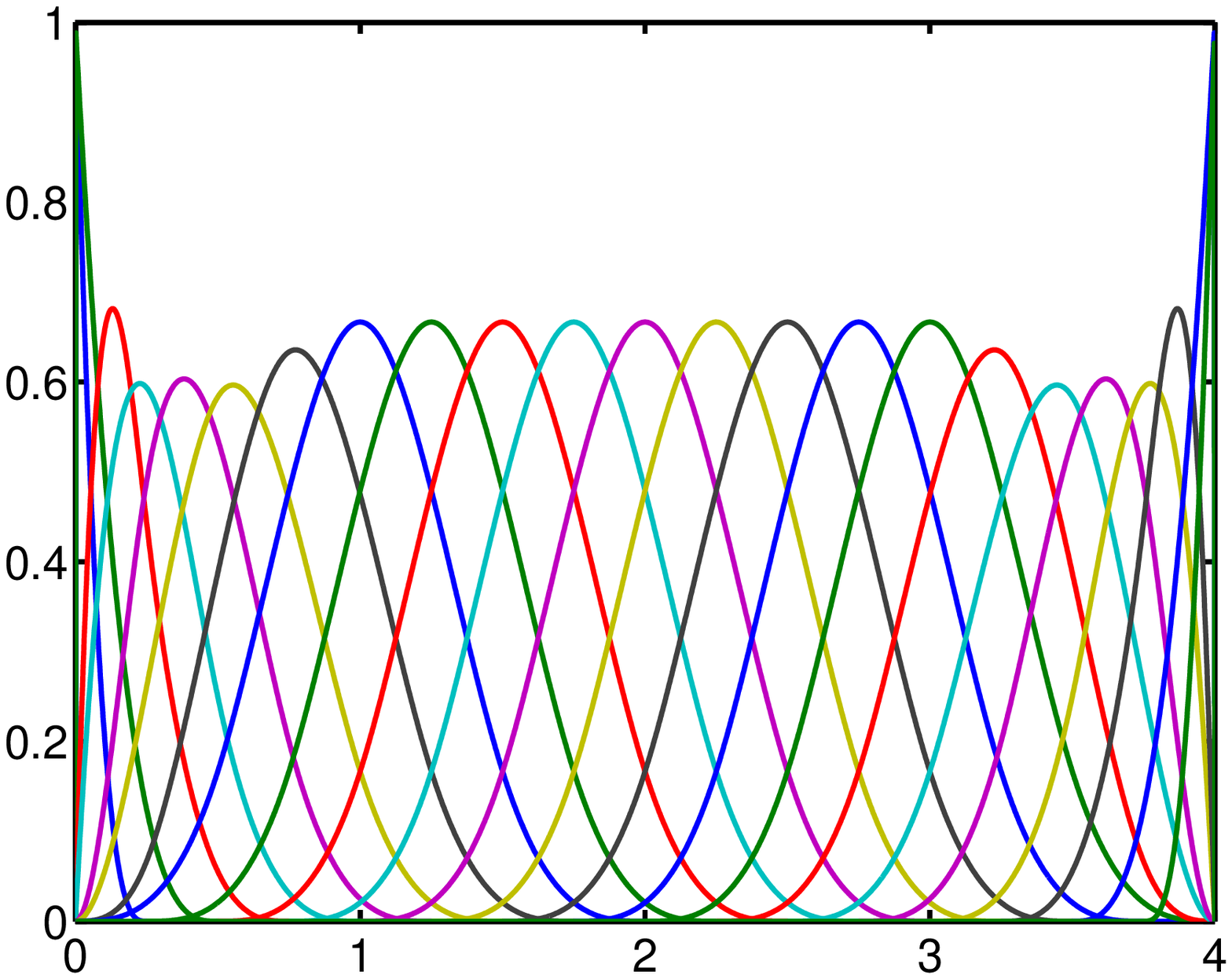}
\end{center}
\caption{Examples of bases (graphs on the left) 
and the corresponding dictionaries consisting of 
functions of double support. The top graphs correspond 
to B-splines of order $m=1$. The middle and 
bottom graphs involve B-splines of order $m=4$ and 
ESEP and EPKB settings, respectively.}
\label{f1}
\end{figure}

For the sake of a simpler notation 
we prefer to build 
dictionaries  by 
considering the extended partition ESEP, 
i.e., $\tDel=\P_b(c-mb,d+mb)$. 
The reason is that,  
since the whole partition is then equidistant, we can
construct the spanning functions  by translation of
 one prototype B-spline of support of length $mb$ into the points
$\{y_i\}_{i=1}^{m+N}= \P_b(c-mb,d)$  
and  the restriction of the functions 
 to the interval $[c,d]$.  
This process is
equivalent to constructing the
boundary functions by truncation 
(see Figure~\ref{f1}), which allows us to use 
a simple notation to label the dictionary functions. 
Shifting indices, for later convenience, we
indicate the basis for $S_m(\Del)$ associated 
to the ESEP $\tDel$ as
\begin{equation}\label{basis1}
\{\phi_k(x)\}_{k\in \Pb}=\{\phi(x-k)\}_{k\in \Pb},
\end{equation}
where $\phi(x)=\frac{1}{b}B\left(\frac xb\right)$. Note that $\supp(\phi)=[0,mb]$.

Since $S_m(\Del)$ is a space of dimension $m+N$, it 
is clear that the dimension can be increased by
decreasing the distance between knots. 
The main contribution of the present effort is to propose
an alternative way of increasing dimension, namely:
maintaining the same distance between 
knots and including functions arising  
by simple translations of  
the basis function for $S_m(\Del)$. 
This is established by the following theorem. 

\begin{theorem}\label{maintheorem}
Let $\Del=\P_b[c,d]$ 
and $\Del'=\P_{b'}[c,d]$
 be such that $\Del\subset\Del'$. Let us denote 
as $\{\phi_k\}_{k \in \Pb}$ and 
$\{\phi'_k\}_{k \in \P_{b'}(c-mb',d)}$ the 
corresponding ESEP B-spline basis 
for $S_m(\Del)$ and $S_m(\Del')$, respectively. 

We construct a {\em dictionary}, $\D_m(\Del,b')$, 
of B-spline
functions on $[c,d]$ 
as 
\begin{align}{\label{rel1}}
\D_m(\Del,b')&=\{\phi_k(x)\}_{k\in \Pb[b']},\\
\intertext{for which it holds that}
\label{main.eq}
\Span\{\D_m(\Del,b')\}&= S_m(\Del').
\end{align}
\end{theorem}
Note that the number of functions in the above 
defined dictionary is equal to the 
cardinality of
$\JDb[b']$, which happens to be
$K=(d-c+mb)/{b'}-1$.
Before advancing the proof of this theorem let us 
assert the following remark:
\begin{remark}\label{remark-S_m}
Setting  $p=b/b'-1$,   for $0 < i \le p$, we have
\begin{align}
S_m(\Del\uplus ib')&=\Span\{\phi_k(x)\}_{k\in \{\P_b[c-mb,d)+ ib'\}}= \Span\{\phi(x-(k+ib'))\}_{k\in\P_b[c-mb,d)},
\label{Sm}\\
\Pb[b']&=\Pb \cup \bigcup_{i=1}^p\{\P_b[c-mb,d)+ ib'\}.
\label{Del2}
\end{align}
Hence, an alternative definition for the dictionary
$\D_m(\Del,b')$, which is equivalent to \eqref{rel1}, 
is the following:
\begin{equation}\label{newdict}
\D_m(\Del,b')=
\{\phi(x-k)\}_{k\in\Pb[b]}
\cup
\bigcup_{i=1}^p \{\phi(x-(k+ib')\}_{k\in\P_b[c-mb,d)}.
\end{equation}
Moreover, from \eqref{Sm} and \eqref {newdict}, 
it follows that 
\begin{equation}\label{dict-union}
\Span\{\D_m(\Del,b')\}=
\bigcup_{i=0}^p S_m(\Del\uplus ib').
\end{equation}
\end{remark}
We are now in a position to start the  proof of 
Theorem~\ref{maintheorem}.
\begin{proof}[Proof of Theorem~\ref{maintheorem}]
We will use \eqref{dict-union}  to prove 
that $S_m(\Del')=\Span\{\D_m(\Del,b')\}$.\\
The inclusion $S_m(\Del')\supset \Span\{\D_m(\Del,b')\}$
follows from \eqref{dict-union} and the 
fact that, from Property~\ref{prop1}, 
$S_m(\Del),S_m(\Del\uplus b'), \dots, S_m(\Del\uplus pb') \subset S_m(\Del')$. 
We can then construct  
scaling-like equations by expressing each  
function $\phi_k\in\D_m(\Del,b'), k \in \JDb[b']$ as a linear 
combination of functions $\phi'_n \in S_m(\Del') , n\in \Pbp$. 
Considering  the support of functions 
$\phi_{k}$  we introduce the scaling equations by 
grouping them into three classes:
\begin{subequations}\label{scal.eq.2}
\begin{align}
  \phi_{k}(x)&=\sum_{n\in\JL}h_{n,k} \phi'_{n}(x)
  &&\text{ for } k\in \JLDb[b'] ,\ &\JL&=\P_{b'}(c-mb',k+m(b-b')], 
\label{first2}\\
  \phi_{k}(x)&=\sum_{n\in\JI}h_{n,k} \phi'_{n}(x)
  &&\text{ for } k\in \JIDb[b'] ,\ &\JI&=\P_{b'}[k,k+m(b-b')],\label{second2}\\
  \phi_{k}(x)&=\sum_{n\in\JR}h_{n,k} \phi'_{n}(x)
  &&\text{ for } k\in \JRDb[b'],\ &\JR&=\P_{b'}[k,d). \label{third2}
\end{align}
\end{subequations}
Equation \eqref{first2} describes the 
{\em left boundary} functions, i.e., all the functions which 
have support of length smaller than $mb$ and contain the point $c$ in
their support. Equation \eqref{second2} 
describes the {\em inner} functions, which are  
functions of  support of length $mb$ arising 
just as translations of a fixed prototype B-spline function.
Finally, equation \eqref{third2} corresponds to 
 the {\em right boundary} functions, which  
have support of length smaller than $mb$ and contain the point $d$ in
their support.

The proof of the inclusion
$S_m(\Del') \subset\Span\{\D_m(\Del,b')\} $
will be achieved by showing that
every $\phi'_{n}(x)\in S_m(\Del'),\ n\in  \Pbp$  can be
expressed as a  linear combination of
functions $\phi_{k}\in\D_m(\Del,b'), k\in \JDb[b']$. 
For such an end we find it convenient to merge  the 
equations \eqref{second2} and \eqref{third2}. Notice that 
 \eqref{second2} 
 and \eqref{third2} can be recast as:
\begin{equation}
\phi_{k}(x)=\sum_{n\in\JR}h_{n,k} \phi'_{n}(x)
 \text{ for }\quad k\in\P_{b'}[c,d) ,\quad \JR=\P_{b'}[k,d). \label{b+c}
\end{equation}
Now, since $\Pbp= \P_{b'}(c-mb',c) \cup \P_{b'}[c,d)$, 
in order to decompose every function  $\phi'_{n},\ n\in  \Pbp$
in terms of the broader functions $\phi_{k}, k\in \JDb[b']$
we make use of \eqref{b+c} if $n \in \P_{b'}[c,d)$ and
\eqref{first2} if  $n \in \P_{b'}(c-mb',c)$. 
The steps leading to the  required 
decompositions are spelled out below 
by focusing on the case  $n \in \P_{b'}[c,d)$.

Let us recall that each 
function $\phi_k,  k\in\P_{b'}[c,d)$ is supported in 
$[k,k+mb]\cap [c,d]$ 
and that 
$\phi_k |_{[k,k+b']}=h_{k,k}\phi'_k|_{[k,k+b']}$.
Hence one can be certain that the coefficient 
$h_{k,k}$ in \eqref{b+c} is non-zero and use 
this equation to express
$\phi'_k(x)$ for $k\in \P_{b'}[c,d)$, as:
\begin{equation}\label{change}
\phi'_{k}(x)=\frac{1}{h_{k,k}}\phi_{k}(x)-
\sum_{n\in\P_{b'}[k+b',d)}\frac{h_{n,k}}{h_{k,k}}\phi'_{n}(x).\end{equation}
The proof that for all $l\in \P_{b'}[c,d)$ each function 
$\phi'_{l}(x)$ is expressible as a linear combination of 
functions $\phi_{k}(x), k\in \P_{b'}[c,d)$ follows  
by subsequent 
evaluations of \eqref{change} at the explicit values 
$k=l,l+b',l+2b',\ldots,d-b'$. Indeed, for $k=l$ we have:
\begin{equation}
\label{inter}
 \phi'_{l}(x)=\frac{1}{h_{l,l}}\phi_{l}(x)-
\sum_{n\in\P_{b'}[l+b',d)}\frac{h_{n,l}}{h_{l,l}}\phi'_{n}(x)
=\frac{1}{h_{l,l}}\phi_{l}(x) - 
\frac{h_{l+b',l}}{h_{l,l}}\phi'_{l+b'}(x) - 
\sum_{n\in\P_{b'}[l+2b',d)}\frac{h_{n,l}}{h_{l,l}}\phi'_{n}(x),
\end{equation}
and the recursive process evolves as follows:
The first step entails to evaluate 
\eqref{change} at $k=l+b'$ 
and introduce the corresponding expression 
for $\phi'_{l+b'}(x)$
in \eqref{inter}. Thus we have:
\begin{equation}\label{decomp}
 \phi'_{l}(x)=\frac{1}{h_{l,l}}\phi_{l}(x)-
\frac{h_{l+b',l}}{h_{l,l}h_{l+b',l+b'}}\phi_{l+b'}(x)+
\sum_{n\in\P_{b'}[l+2b',d)}\left(
\frac{h_{l+b',l}h_{n,l+b'}}{h_{l,l}h_{l+b',l+b'}}
-\frac{h_{n,l}}{h_{l,l}}\right)\phi'_n(x).
\end{equation}
The next step consists 
of evaluating \eqref{change} at $k=l+2b'$ 
and introducing the corresponding expression for
$\phi'_{l+2b'}(x)$
in \eqref{decomp}. The process continues repeating 
equivalents steps, i.e., at step $s$ say, we 
evaluate \eqref{change} at $k=l+sb'$ and introduce 
the equation in \eqref{decomp}. 
Finally, at step $(d-b'-l)/b'$, for $k=d-b'$ we obtain: 
\begin{equation}
\label{last}
\phi'_{d-b'}(x)=
\frac{1}{h_{d-b',d-b'}}\phi_{d-b'}(x).
\end{equation}
When $\phi'_{d-b'}(x)$ as given in \eqref{last}
is introduced in \eqref{decomp}, the right hand 
side of this equation turns out to be a linear 
combinations of functions $\phi_{l}(x), \phi_{l+b'}(x),
\phi_{l+2b'}(x), \ldots, \phi_{d-b'}$. We have thereby 
proved that   $\phi'_{l}(x)\in 
\Span\{\D_m(\Del,b')\}, l\in \P_{b'}[c,d)$. 

The proof concerning functions 
$\phi'_{n}(x),\ n\in \P_{b'}(c-mb',c)$ parallels the proof 
given above, but the following considerations are in order: 
To obtain the corresponding 
equation for $\phi'_{n}(x)$ we use 
\eqref{first2}
and the fact that  
$h_{n,k}\ne 0$ for $k=n-m(b-b')$. The required 
decomposition arises by 
evaluating such an equation at values of $k$ 
in decreasing order
ranging from $k=n-m(b-b')$ to $k=c-mb+b'$.  
\end{proof}
\begin{corollary}
The dictionary
$\D_m(\Del,b')=\{\phi_k(x)\}_{k\in \JDb[b']}$ is a
frame for $S_m(\Del')$, i.e.,
for all $f \in S_m(\Del')$ there exist two constants
$0<A\le B$ such that
\begin{equation}
 A || f||^2 \le \sum_{k\in \JDb[b']}|\langle f,\phi_k \rangle|^2\le B||f||^2.
\end{equation}
\end{corollary}
\begin{proof} The proof is a direct consequence of
 Theorem~\ref{maintheorem} since the dictionary $\D_m(\Del,b')$ is a
finite dimension set of functions of finite norm 
and, as such, a frame for its span. The upper bound 
is a consequence of Schwartz inequality and the non-zero 
lower bound is ensured by the fact that for 
$f \in S_m(\Del')$ it is true that 
 $\langle f,\phi_k \rangle \ne 0$ at least for 
 one $k$ in $\JDb[b']$. 
\end{proof}
\begin{remark} It is appropriate to point out that, 
although the proof of Theorem~\ref{maintheorem} was given 
for dictionaries constructed within an 
ESEP setting, the result can be extended to more  
general extended  partition settings. 
The corresponding proof is equivalent to the one given here, but it involves a less 
handy notation. 
\end{remark}
Figure \ref{f1} shows some examples of 
B-spline dictionaries. 
The top graph on the left depicts a B-spline basis of 
order $m=1$.  
These are piecewise constant functions and, since
there are not boundary functions, the ESEP and EPKB basis 
are equivalent. 
The corresponding dictionary 
of functions of double support is depicted in 
the top right graph of the same figure. As  
emphasized by the thicker lines,
the dictionary contains two boundary functions 
of support equal to the basis function of the 
left graph. 
The middle left graph of Figure~\ref{f1} 
shows the B-spline basis of order 
$m=4$ corresponding to the 
ESEP extended partition. The middle right graph 
shows the dictionary, for the same space, consisting 
of functions of twice as much  support as the 
corresponding basis functions. 
The graphs at the bottom have the same description as 
the middle graphs, but  correspond 
to the EPKB  extended partition.

\section{Application to recursive signal approximation}
\label{sec:application}

We present here two examples to illustrate
the relevance of the proposed dictionaries
to a typical problem of signal representation:
the problem of achieving a sparse
atomic decomposition approximating a given signal.

\begin{figure}[!ht]
\begin{center}
\includegraphics[width=6cm]{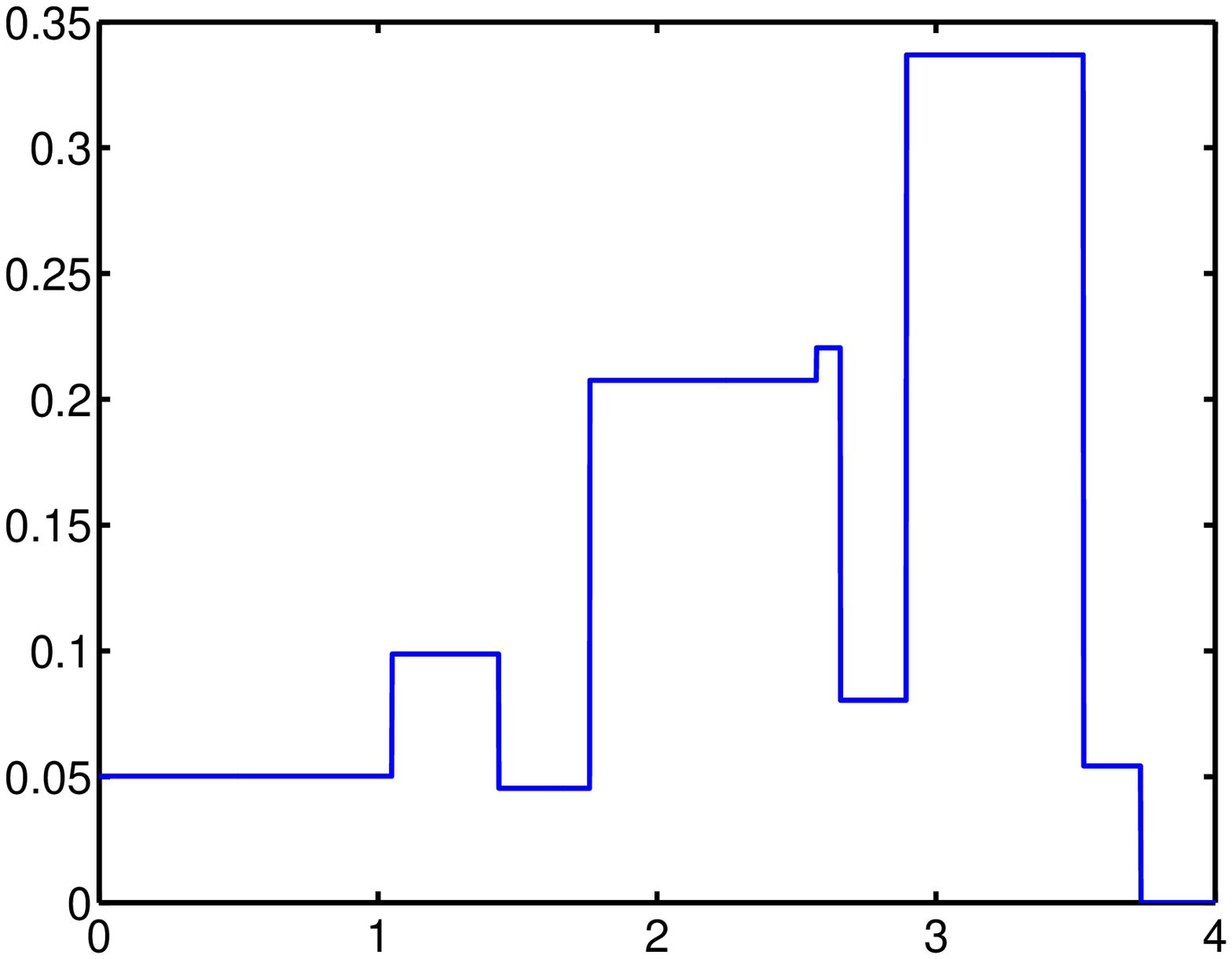}\qquad
\includegraphics[width=6cm]{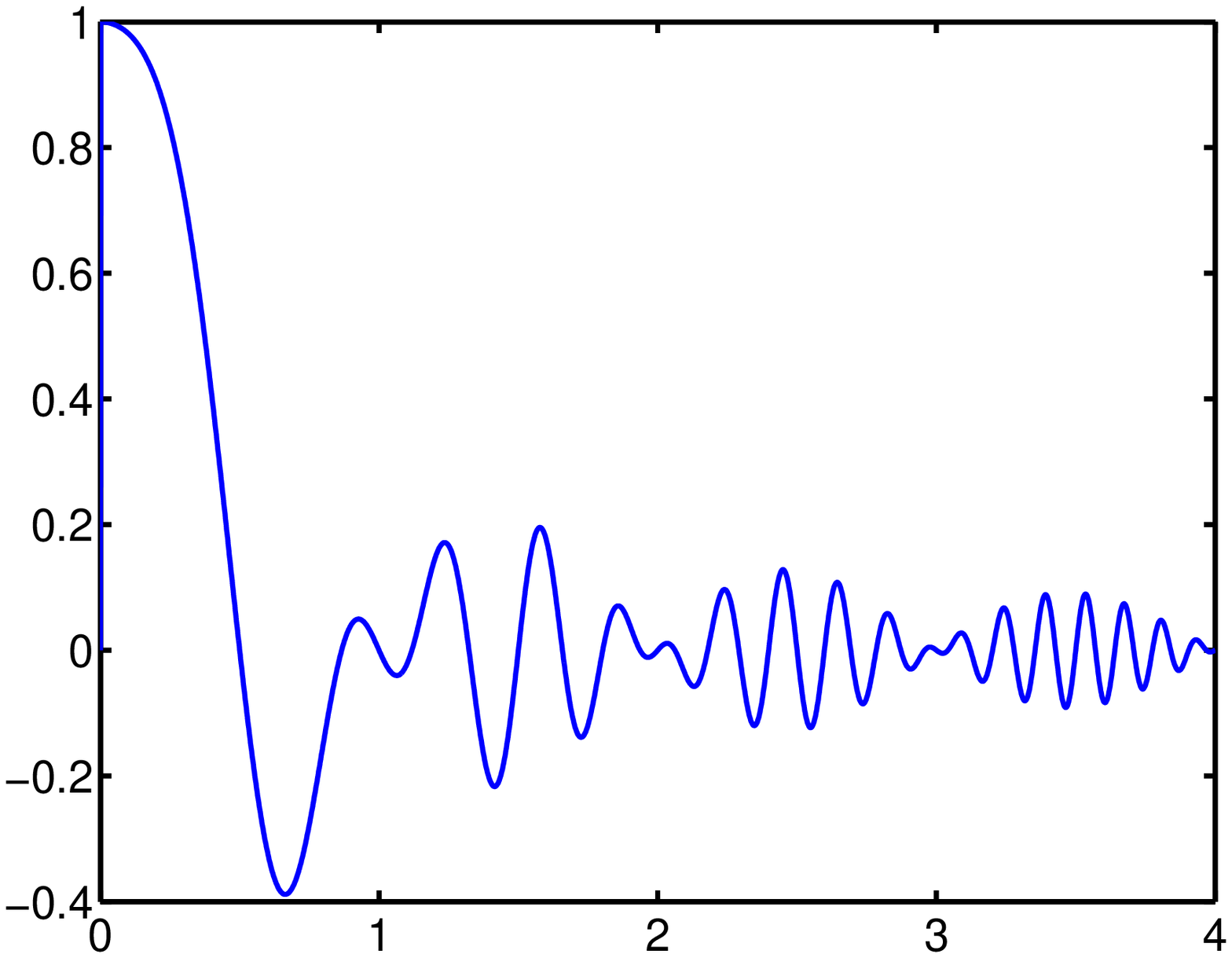}
\end{center}
\caption{Randomly generated blocky signal (left graph). 
Modulated chirp (right graph).}
\label{f2}
\end{figure}

As a first example we consider the
randomly generated blocky signal
on the interval $[0,4]$ depicted
in the left graph of Figure~\ref{f2}.  Such a
signal has an acceptable approximation
in the subspace $S_1(\Del') $,  with    
 distance between knots $b'=2^{-8}$, spanned by a B-spline basis
of order one. 
The graph indicating the signal approximation 
 coincides with 
 the one of Figure~\ref{f2}. For such an approximation 
 we need to use
 almost all the basis functions  (956 out of 1024).
 This is of course an expected result, since the
 support of the basis functions for $ S_1(\Del') $ is
 considerably small ($|\supp|=2^{-8}$) in comparison to the
 length of the blocks composing the signal.
 It is then convenient to construct the
 identical approximation using 
 dictionary functions of much larger 
 support ($|\supp|=1$) spanning the same space  $ S_1(\Del') $.
 In this case the approximation is obtained  
 by means of only $M=23$ functions, 
 chosen through a
 forward and backward optimized orthogonal matching
 pursuit strategy \cite{boomp,oomp}, from a 
 dictionary of $1279$ atoms.
 For further comparison
 we have used the Haar wavelet representation
 for the same space, 
 which contains wavelets 
 of varied support and scaling 
 functions of the same support as the 
 dictionary functions.
 Since Haar wavelets are orthonormal, to obtain the
 desired signal approximation we need just to retain
 $M=50$ basis functions corresponding to the 50 
 coefficients of largest absolute value.

The second example involves the piece of 
 modulated chirp signal plotted
 in right graph of Figure~\ref{f2}.
This signal has
an acceptable representation 
in the subspace $S_4(\Del')$,  
with $b'=2^{-5}$, 
spanned by the B-splines  basis of order
four.
Since in this case the basis 
is non-orthogonal, in order to
 approximate the signal by a subset of basis functions
  we applied the same pursuit strategy as in the
  previous example. An approximation coinciding with the
  graph of Figure~\ref{f2} is obtained
 with $M=101$ basis functions. Using a dictionary of functions having
 four times larger support (i.e., distance $2^{-3}$ between knots), and the same matching
 pursuit strategy for selecting functions, we need 
 $M=58$ dictionary functions. In this case, the selection from the 
 wavelet representation of the subspace needs to be carried out also 
 through the matching pursuit approach. The required number of 
 wavelet functions is $61$.

It is interesting to notice that in both examples we
have achieved representations which in terms of
sparseness are comparable (superior in the first case)
to wavelet basis representation. This is a surprising
result, since wavelet basis
are composed by functions of different support and our
dictionaries by functions of
 fixed support. These outcomes are
 certainly worth to be investigated further.

MATLAB codes for generating  the proposed 
dictionaries and implementing pursuit strategies are available at 
\cite{webpage}.
 
\section{Conclusions}
\label{sec:conclusion}
An interesting feature of B-spline functions
ha s been discussed. We have shown  that 
a dictionary for a cardinal spline space on a 
compact interval can be constructed by translating a 
prototype B-spline function into the knots of the 
corresponding space. This property allows to 
span a cardinal spline space by using B-spline 
functions of larger support than the one corresponding 
to the basis functions for the same space. 

As an example of application
of the proposed dictionaries, 
two signals of different
nature have been represented by selecting atoms
from dictionaries of B-splines.
The results illustrate the possible relevance of
B-spline dictionaries
to problems requiring sparse representation.

Finally, we would like to 
point out some lines of follow 
up work that we believe to be interesting: 
The possibility of spanning a fixed space by 
B-spline dictionaries, each of which consists of 
 functions of different support, 
 arises the question as to how to
choose in a effective manner the dictionary of
B-splines of `optimal support' for
representing a given signal. 
Another matter that
appears definitely worth looking at
is the possibility of extending the proposed 
construction in order to generate suitable 
subspaces by translating wavelets generated 
by B-splines.

\end{document}